\documentclass[10pt,leqno]{amsart}
\usepackage{amscd}
\usepackage{amssymb}
\usepackage{amsfonts}
\usepackage{latexsym}
\usepackage{verbatim}

\theoremstyle{plain}
\newtheorem{theorem}{Theorem}[section]
\newtheorem*{theorem*}{Theorem}
\newtheorem{definition}[theorem]{Definition}
\newtheorem{lemma}[theorem]{Lemma}
\newtheorem{prop}[theorem]{Proposition}
\newtheorem{cor}[theorem]{Corollary}
\newtheorem{rem}[theorem]{Remark}
\newtheorem{ex}[theorem]{Example}
\begin{document}
\title{On Special Calibrated Almost Complex Structures and Moduli Space}
\author{Adriano Tomassini and Luigi Vezzoni}
\date{\today}
\address{Dipartimento di Matematica\\ Universit\`a di Parma\\ Viale G.
  P. Usberti 53/A\\
43100 Parma\\ Italy}
\email{adriano.tomassini@unipr.it}
\address{Dipartimento di Matematica \\ Universit\`a di Torino\\
Via Carlo Alberto 10 \\
10123 Torino\\ Italy} \email{luigi.vezzoni@unito.it}

\thanks{This work was supported by the Project M.I.U.R.
``Geometric Properties of Real and Complex Manifolds'' and by G.N.S.A.G.A.
of I.N.d.A.M.}
\begin{abstract}
An \emph{$\omega$-admissible almost complex structure} on a
$2n$-dimensional symplectic manifold $(M,\omega)$ is a
$\omega$-calibrated almost complex structure $J$ admitting a
nowhere vanishing $\overline{\partial}_J$-closed $(n,0)$-form
$\psi$. After giving some examples we consider the moduli space of
admissible almost complex structures and we study infinitesimal
deformations. As special case, we write down explicit computations for
the complex torus.
\end{abstract}
\maketitle
\newcommand\C{{\mathbb C}}
\newcommand\R{{\mathbb R}}
\newcommand\Z{{\mathbb Z}}
\newcommand\T{{\mathbb T}}
\newcommand{\de}[2]{\frac{\partial #1}{\partial #2}}
\newcommand\wJ{{\widetilde{J}}}
\newcommand{\ov}[1]{\overline{ #1}}
\newcommand{\Tk}{\mathcal{T}_{\omega}}
\newcommand{\op}{\overline{\partial}}
\section{Introduction}
The interplay between symplectic and complex structures has been
studied quite extensively in the last years (see e.g. \cite{AD},
\cite{AL}, \cite{G}, \cite{GHJ}, \cite{J1} and the references therein). It is well known
that on an a symplectic manifold $(M,\omega)$ there exist
many almost complex structures $J$ such that
\begin{equation}\label{tamed}
\omega[x](v,w)=\omega[x](J_xv,J_xw)\,,\quad
\omega[x](u,J_xu)>0 \,,
\end{equation}
for any $x\in M\,,$ $v,w,u\in T_xM, \,u\neq0$ (see e.g. \cite{AL},
\cite{G}). An almost complex structure $J$ on a
symplectic manifold $(M,\omega)$ is said to be $\omega$-{\em
calibrated} if it satisfies condition \eqref{tamed}. Furthermore
the existence of special holomorphic structures on a symplectic
manifold (e.g. K\"ahler structures, Calabi-Yau structures) imposes
strong conditions on the topology of the manifold. Hence, it is
natural to consider the non-integrable case in order to have more
flexible structures. In this context, the notion of symplectic
Calabi-Yau manifold has been considered in \cite{dB2}, \cite{dBT}
and it appears as a natural generalization of Calabi-Yau
manifold\footnote{In \cite{dBT} symplectic Calabi-Yau manifolds were called
\emph{generalized Calabi-Yau manifolds}.
Here we change the terminology to avoid confusion with the case considered by Hitchin in \cite{H}}.
This generalization is  different from that one given by Hitchin
in \cite{H} in the context of \emph{generalized geometry}.

Namely, a \emph{symplectic Calabi-Yau manifold} is a symplectic manifold
$(M,\omega)$ endowed with an $\omega$-calibrated almost complex structure
$J$ and a complex volume form $\psi$ covariantly constant with
respect to the Chern connection of $(M,J,\omega)$. This is
equivalent to have an almost K\"ahler manifold
endowed with a $\op_J$-closed complex volume form. In dimension $6$ this definition can be improved by
requiring that the real part of the complex volume form be closed. Such structures are called
\emph{symplectic Half-flat}.\\

In the present paper we study the moduli space of $\omega$-calibrated almost complex structures admitting a
symplectic Calabi-Yau structure. Such almost complex structures will be called $\omega$-\emph{admissible}.

First of all, we
give
an example of a non-admissible almost complex structure on a $6$-dimensional compact symplectic manifold.
This manifold provides an example of an almost K\"ahler manifold whose Chern connection has holonomy
non-contained in $\mbox{SU}(3)$. Then,
we study the infinitesimal deformations of the space of
$\omega$-admissible almost
complex structures and we compute the virtual tangent space to the moduli space.
In the last part we apply our results
to the torus showing that the standard complex structures is not rigid.\\

In section 2 we recall some preliminary results on complex manifolds and Hermitian geometry and give
the basic definitions.
In section 3 we give the example described above and an example of
a compact almost complex $6$-manifold which admits symplectic Calabi-Yau structures,
but has no symplectic half-flat structure.

In section 4 we study the infinitesimal deformations of the space of
$\omega$-admissible almost complex structures $\mathcal{AC}_\omega(M)$. We introduce the deformation form $\theta_L$
(see subsection
\ref{santideformationform} for the precise definition), which is a $1$-form depending on the choice of an endomorphism
of $TM$ and on a complex $\op$-closed volume form, and we prove the following
\begin{theorem*}
Let $(M,\omega)$ be a symplectic manifold. Fix $J
\in\mathcal{AC}_\omega(M)$ and consider a smooth curve $J_t$ in
$\mathcal{AC}_\omega(M)$ close to $J$ and satisfying $J_0=0$. Then
the derivative $\dot{J}_0$ of $J_t$ at $0$ is given by
$$
\dot{J}_0=2\,JL,
$$
where
$L\in\mbox{\emph{End}}(TM)$ satisfies  the following conditions
$$
L=\,^t\!L\,,\quad LJ=-JL
$$
and for any nowhere vanishing $\op_J$-closed
$\psi\in\Lambda^{n,0}_J(M)$ the deformation form $\theta_L(\psi)$ is $\op_J$-exact.
\end{theorem*}
A key tool to prove this theorem is proposition
\ref{santiprop}, which describes the relationship between the $\op$-operators of
two close almost complex structures (see section 4).\\

In the compact case the previous theorem allows to define the virtual tangent space to the
moduli space of admissible almost complex structures
$$
\mathfrak{M}(\mathcal{AC}_\omega(M))=\mathcal{AC}_{\omega}(M)/\mbox{Sp}_\omega(M)\,.
$$
In particular we can introduce the concepts of \emph{unobstructed} and \emph{rigid}
$\omega$-admissible almost complex structures. See proposition \ref{centrale} and definition \ref{rigido}.

In section 5 we apply the results of section 4 to
to the complex torus, showing that its standard complex
structure is not rigid.\\ \newline
This paper has originated from a series of seminars given in
Florence. We would like to thank Paolo de Bartolomeis for useful comments and remarks. We also thank Paul
Gauduchon for his precious help.

\section{Preliminaries}
Let $(M,J)$ be an almost complex manifold manifold. Then the bundle
of complex valued $r$-forms $\Lambda^r_{\C}(M)$ decomposes as
$$
\Lambda^{r}_{\C}(M)=\bigoplus_{p+q=r}\Lambda^{p,q}_J(M)\,,
$$
where $\Lambda^{p,q}_J(M)$ is the bundle of $(p,q)$-forms on
$(M,J)$. According with the above decomposition, the exterior
derivative
$$
d\colon \Lambda^r_{\C}(M)\to \Lambda^{r+1}_{\C}(M)
$$
splits as
$$
d\colon\Lambda^{p,q}_J (M)\to
\Lambda^{p+2,q-1}_J(M)\oplus\Lambda^{p+1,q}_J(M)\oplus\Lambda^{p,q+1}_J(M)\oplus\Lambda^{p-1,q+2}_J
(M)\,,
$$
$$
d=A_J+\partial_J+\op_J+\overline{A}_J\,.
$$
Let $\omega$ be an almost symplectic structure on $M$, i.e. $\omega$ is a non-degenerate $2$-form on $M$.
An almost
complex structure $J$ is said to be $\omega$-\emph{tamed} if
$$
\omega[x](v,J_xv)>0
$$
for any $x\in M$, $v\in T_xM$, $v\neq 0$.

If $J$ satisfies the extra condition
$$
\omega[x](J_xv,J_xw)=\omega[x](v,w)\,,
$$
for any $x\in M$, $v,w\in T_xM$, then $J$ is said to be
$\omega$-\emph{calibrated}. In this case
$$
g_J[x](v,w):=\omega[x](v,J_xw)\,,
$$
defines an almost Hermitian metric on $M$ (i.e. a $J$-invariant Riemannian
metric).  We denote by $\mathcal{T}_{\omega}(M)$ and by
$\mathcal{C}_\omega(M)$ the space of $\omega$-tamed and
$\omega$-calibrated almost complex structures on $M$ respectively.
It is well known that $\mathcal{T}_{\omega}(M)$ is a contractible
space (see e.g. \cite{AL}). In particular, the first Chern class of
$(M,J)$ does not depend on the choice of $J\in
\mathcal{T}_{\omega}(M)$. Hence it is well defined the first Chern
class of $(M,\omega)$.

By definition, an almost complex structure $J$ is said to be \emph{integrable} (or a \emph{complex structure}) if the
Nijenhuis tensor
$$
N_J(X,Y)=[JX,JY]-J[JX,Y]-J[X,JY]-[X,Y]
$$
vanishes.\\
\newline
In the sequel we will use the following fact (see e.g. \cite{HE}): \vskip0.2truecm
\noindent\emph{let $(M,J)$ be a compact almost complex manifold
and $f\colon M\to \C$ be a $J$-holomorphic function, then $f$ is
constant.} \vskip0.2truecm \noindent
Now we are going to recall the definition of {\em symplectic
Calabi-Yau structure}. Let $(M,\omega,J)$ be an almost K\"ahler
manifold and let $\nabla^{LC}$ be the Levi Civita connection of
the metric $g_J$ associated with $(\omega,J)$. Then the Chern
connection is defined as
$$
\nabla=\nabla^{LC}-\frac{1}{2}J\nabla^{LC} J
$$
and it satisfies the following properties
$$
\nabla g=0\,,\quad \nabla J=0\,,\quad T^{\nabla}=\frac{1}{2}\,N_J\,,
$$
$T^{\nabla}$ being the torsion of $\nabla$.
\begin{rem}
\label{n=op}
\emph{We recall that $\nabla^{0,1}=\op_J$ (see e.g. \cite{Gau}), where $\nabla^{0,1}$ denotes the
$(0,1)$-component of $\nabla$.}
\end{rem}
We have the following definition (see \cite{dBT})
\begin{definition}
A \emph{symplectic Calabi-Yau} manifold consists of
$(M,\omega,J,\psi)$, where $(M,\omega)$ is a $2n$-dimensional
symplectic manifold, $J$ is an $\omega$-calibrated almost complex
structure on $M$ and $\psi$ is a nowhere vanishing $(n,0)$-form on
$M$ satisfying
$$
\nabla\psi=0\,,
$$
where $\nabla$ is the Chern connection of  $(\omega,J)$.
\end{definition}
\begin{rem}
\emph{Note that a Calabi-Yau manifold is in particular a
symplectic
  Calabi-Yau manifold. Indeed in this case the Chern connection is the
  Levi Civita one, since $J$ is integrable.}
\end{rem}
In the 6-dimensional case we can improve the previous definition
by requiring that the real part of $\psi$ is $d$-closed; namely
\begin{definition}\label{SGCY}
A \emph{symplectic half-flat manifold} is the datum of
$(M,\omega,J,\psi)$, where $(M,\omega)$ is a $6$-dimensional
symplectic manifold, $J$ is an $\omega$-calibrated almost complex
structure on $M$, $\psi$ is a nowhere vanishing $(3,0)$-form such
that
$$
\nabla\psi=0\,,\quad
\psi\wedge\overline{\psi}=-\frac{4}{3}i\,\omega^3 \,,\quad d\Re
\mathfrak{e}\,\psi=0\,.
$$
\end{definition}
These structures are just the intersection between symplectic and
half-flat ones. The latter have been introduced and studied by
Hitchin and Chiossi-Salamon (see \cite{H1} and \cite{CS}). In this
situation is possible to perform (special) Lagrangian geometry by
considering Lagrangian submanifolds calibrated by
$\Re\mathfrak{e}\,\psi$ (see \cite{dB2}).
\begin{rem}
\emph{Note that the condition $\nabla\psi=0$ is redundant since it
  can be showed that, given a nowhere vanishing
  $\psi\in\Lambda^{3,0}_J(M)$,  then the following facts are equivalent  (see \cite{dBT})
$$
\begin{cases}
&\psi\wedge\ov{\psi}=i\lambda\,\omega^3\\
&d\Re \mathfrak{e}\,\psi=0
\end{cases}
\iff
\begin{cases}
&\nabla\psi=0\\
&(A_{J}+\ov{A}_J)\psi=0\,,
\end{cases}
$$
where $\lambda\in\R$, $\lambda\neq 0$ and
$A_J\colon\Lambda^{p,q}_J(M)\to \Lambda^{p+2,q-1}_J(M)$ denotes the $(p+2,q-1)$-component of the exterior derivative of
a $(p,q)$-form.
}
\end{rem}
We introduce now the definition of {\em admissible almost complex
structures}. Namely, if $(M,\omega)$ is a symplectic manifold, an
$\omega$-calibrated almost complex structure $J$ on $M$ will be
called \emph{admissible} if there exists
$\psi\in\Lambda_J^{n,0}(M)$ such that the triple $(\omega,J,\psi)$
is a symplectic Calabi-Yau structure on $M$. More precisely we
state the following
\begin{definition}
Let $(M,\omega)$ be a symplectic manifold. An almost complex
structure $J$ on $M$ is said to be $\omega$-\emph{admissible} if
\begin{enumerate}
\item[1.] $J$ is $\omega$-calibrated; \item[2.] there exists a
nowhere vanishing $\psi\in\Lambda^{n,0}_J(M)$ such that
$\op_J\psi=0$.
\end{enumerate}
\end{definition}
In the next sections we will study some properties of admissible
almost complex structures and we will introduce the Moduli space
of such structures.
\section{Examples}
We start by giving an example of an almost complex structure which
admits a symplectic Calabi-Yau structure, but it has no symplectic
half-flat structures.
\begin{ex}\label{ex}
\emph{
Let
$$
G=\left\{
\left(
\begin{array}{ccc}
1  &z_1  &z_2\\
0  &1    &z_3\\
0  &0    &1
\end{array}
\right)\,\,|\,\,z_1,z_2,z_3\in\C\right\}
$$
be the complex Heisenberg group and let $\Gamma\subset G$ be the
subgroup with integral entries. Then $M=G/\Gamma$ is the \emph{Iwasawa
manifold}. It
is known that $M$ is symplectic, but it has no K\"ahler structures
(see \cite{CFG}).\\
Let $z_r=x_r+ix_{r+3}$, $r=1,2,3$, and set
\begin{equation*}
\begin{aligned}
&\alpha_1=dx_1\,,\;\;\alpha_2=dx_3-x_1dx_2+x_4dx_5\,,\;\;\alpha_3=dx_5,\\
&\alpha_4=dx_4\,,\;\;\alpha_5=dx_2\,,\;\;\alpha_6=dx_6-x_4dx_2-x_1dx_5\,,
\end{aligned}
\end{equation*}
then $\{\alpha_1,\dots,\alpha_6\}$ are invariant $1$-forms on $G$,
so that $\{\alpha_1,\dots,\alpha_6\}$ is a global coframe on $M$.
We immediately get
\begin{equation*}
\begin{cases}
d\alpha_1=d\alpha_3=d\alpha_4=d\alpha_5=0\\
d\alpha_2=-\alpha_1\wedge\alpha_5-\alpha_3\wedge\alpha_4\\
d\alpha_6=-\alpha_4\wedge\alpha_5-\alpha_1\wedge\alpha_3\,.
\end{cases}
\end{equation*}
Let $\{\xi_1,\dots, \xi_6\}$ be the dual frame of
$\{\alpha_1,\dots,\alpha_6\}$; then
$$
\begin{cases}
\begin{aligned}
&J(\xi_r)=\xi_{r+3}\,&&r=1,2,3\\
&J(\xi_{3+r})=-\xi_{r}\,&&r=1,2,3\,,
\end{aligned}
\end{cases}
$$
defines an almost complex structure on $M$ calibrated by the symplectic form
$$
\omega=\alpha_{1}\wedge\alpha_{4}+\alpha_{2}\wedge\alpha_{5}+\alpha_{3}\wedge\alpha_{6}\,.
$$
Let
$\psi=(\alpha_1+i\,\alpha_{4})\wedge(\alpha_{2}+i\,\alpha_{5})\wedge(\alpha_{3}+i\,\alpha_6)$,
then a direct computation gives
$$
\begin{cases}
&\psi\wedge\ov{\psi}=-\frac{4}{3}i\,\omega^3\\
&\op_J\psi=0\,.
\end{cases}
$$
Hence the conditions above and remark \ref{n=op} imply
$$
\nabla\psi=0\,.
$$
Hence $(\omega,J,\psi)$ is a symplectic Calabi-Yau structure on
$M$.}

\emph{
Now we prove that there are no nowhere vanishing (3,0)-forms $\eta$ on $M$ such that
$$
d\Re \mathfrak{e}\,\eta=0\,.
$$
In particular $(M,J)$ does not admit any symplectic half-flat
structure.\\
In order to show this let $\eta\in\Lambda^{3,0}_J(M)$; then there
exists $f\in C^{\infty}(M,\C)$ such that
$$
\eta=f\,\psi\,.
$$
Let $f=u+iv$ and set
$$
du=\sum_{i=1}^nu_i\,\alpha_i\,,\quad dv=\sum_{i=1}^n v_i\,\alpha_i\,.
$$
A direct computation shows that
$$
\begin{aligned}
d\Re\mathfrak{e}\,\eta=&(u_6+v_3)\,\alpha_{3456}+(v_2+u_5)\,\alpha_{2456}
+(-u_{1}+v_4)\,\alpha_{1345}+\\
&+(-u_6-v_3)\,\alpha_{1236}
(-u_3+v_6)\,\alpha_{2346}+(-u_5-v_2)\,\alpha_{1235}+\\
&+(u_1-v_4)\,\alpha_{1246}
+(u_3-v_6)\,\alpha_{1356}+v\,\alpha_{1245}+\\
&-v\,\alpha_{1346}
+(u+u_4+v_1)\,\alpha_{1456}
+(u-u_4-v_1)\,\alpha_{1234}+\\
&+(-u_2+v_5)\,\alpha_{2345}
+(u_2-v_5)\,\alpha_{1256}\,,
\end{aligned}
$$
where $\alpha_{ijhk}=\alpha_i\wedge\alpha_j\wedge\alpha_h\wedge\alpha_k$.
Hence $d\Re \mathfrak{e}\,\eta=0$ if and only if $u=v=0$.}$\quad\Box$
\end{ex}
The next nilmanifold provides an example of a compact almost
complex manifold with vanishing first Chern class and such that
there are no Hermitian metrics whose Chern connection has holonomy
contained in SU$(3)$. This is in contrast with the integrable
case, in view of Calabi-Yau theorem. For other examples involving the Bismut connection see \cite{FG}.
\begin{ex}\emph{
Let
$$
H(3)=\left\{ \left(
\begin{array}{ccc}
1  &x  &y\\
0  &1    &t\\
0  &0    &1
\end{array}
\right)\,\,:\,\,x,y,t\in\R\right\}
$$
be the $3$-dimensional Heisenberg group and let $\Gamma\subset
H(3)$ be the cocompact lattice of matrices with integral entries.
Then $X=S^1\times H(3)/\Gamma$ is called the
\emph{Kodaira-Thurston} manifold.\\
Let $M=\T^2\times X$, where $\T^2$ is the $2$-dimensional standard
torus. Then $M$ can be viewed as $\R^{6}/\!\sim$, where
$$
[(x_1,x_2,x_3,x_4,x_5,x_6)]=[(x_1+m_1,x_2+m_2,x_3+m_3,x_4+m_4,x_5+m_5,x_6+m_4x_5+m_6)]\,,
$$
for any  $(m_1,m_2,m_3,m_4,m_5,m_6)\in \Z^6$. The 1-forms
\begin{equation*}
\begin{aligned}
&\alpha_1=dx_1\,,\;\;\alpha_2=dx_2\,,\;\;\alpha_3=dx_3\,,\\
&\alpha_4=dx_4\,,\;\;\alpha_5=dx_5\,,\;\;\alpha_6=dx_6-x_4dx_5\,,
\end{aligned}
\end{equation*}
define a global coframe on $M$. We have
$$
\begin{aligned}
&d\alpha_i=0\,,\,\,\mbox{for }i=1,\dots,5\,,\\
&d\alpha_6=-\alpha_4\wedge\alpha_5\,.
\end{aligned}
$$
The 2-form
$$
\omega=\alpha_{1}\wedge\alpha_{2}+\alpha_{3}\wedge\alpha_{4}+\alpha_{5}\wedge\alpha_{6}\,,
$$
is a symplectic structure on $M$. Let $J$ be the almost complex structure
defined on the dual frame of $\{\alpha_1,\dots,\alpha_6\}$ by the
relations
$$
\begin{aligned}
& J(X_1):=X_2,\quad &&J(X_3):=X_4,\quad
&&J(X_5):=X_6,\\
& J(X_2):=-X_1,\quad &&J(X_4):=-X_3,\quad
&&J(X_6):=-X_5\,.
\end{aligned}
$$
$J$ is an $\omega$-calibrated almost complex structure on $M$. The
form
$$
\psi:=(\alpha_1+i\alpha_2)\wedge(\alpha_3+i\alpha_4)\wedge(\alpha_5+i\alpha_6)
$$
is a nowhere vanishing section of $\Lambda^{3,0}_J(M)$. We have
$$
\overline{\partial}_J\psi=-(\alpha_3-i\alpha_4)\wedge\psi\,.
$$
Now we prove that there are no nowhere vanishing (3,0)-forms $\eta$
on $M$ such that $\op_J\eta=0$. Set
$$
Z_{j}=\frac{1}{2}(X_j-iJX_j)\,,\quad j=1,2,3
$$
and
$$
\zeta_{j}=\alpha_j+iJ\alpha_j\,,\quad j=1,2,3\,.
$$
Any $\eta\in \Lambda^{3,0}_J(M)$ is a multiple of $\psi$, i.e.
there exists a function $f=u+iv\in C^{\infty}(M,\C)$ such that
$\psi=f\,\eta$. Then
$$
\op_J\eta=\op_J(f\psi)=\op_Jf\wedge\psi+f\op_J\psi
=(\sum_{j=1}^3\ov{Z}_{j}(f)\,\ov{\zeta}_{j}-f\,\ov{\zeta}_2)\wedge\psi\,.
$$
Therefore $\op_J\eta=0$ if and only if the following systems of PDE's
are satisfied:
\\[5pt]
\begin{enumerate}
\item[a.]
$
\begin{cases}
\partial_{x_1}u-\partial_{x_2}v=0\\
\partial_{x_2}u+\partial_{x_1}v=0\,,
\end{cases}
$\\
[5pt]
\item[b.]
$
\begin{cases}
\partial_{x_3}u-\partial_{x_4}v-u=0\\
\partial_{x_4}u+\partial_{x_3}v-v=0\,,
\end{cases}
$
\\[5pt]
\item[c.]
$
\begin{cases}
\partial_{x_6}v-\partial_{x_5}u-x_4\partial_{x_6}u=0\\
\partial_{x_6}u+\partial_{x_5}v+x_4\partial_{x_6}v=0\,,
\end{cases}
$
\\[5pt]
\end{enumerate}
where the unknowns $u,v$ are functions on $\R^6$ satisfying
$$
\begin{aligned}
&u(x)=u(x_1+m_1,x_2+m_2,x_3+m_3,x_4+m_4,x_5+m_5,x_6+m_4x_5+m_6)\,,\\
&v(x)=v(x_1+m_1,x_2+m_2,x_3+m_3,x_4+m_4,x_5+m_5,x_6+m_4x_5+m_6)
\end{aligned}
$$
for any $m_1,\dots,m_6\in\Z$.}

\emph{Equations a. imply that $f=f(x_3,x_4,x_5,x_6)$. Since $f$ is a
function on $M$, then $f$ is $\Z-$periodic in $x_3,x_5,x_6$. Therefore
we can take the Fourier expansion of $u$ and $v$.
Set
$$
\begin{aligned}
&u(x_3,x_4,x_5,x_6)=\sum_{N\in\Z^3} u_{N}(x_4)\, e^{2\pi i(n_3x_3+n_5x_5+n_6x_6)}\,,\\
&v(x_3,x_4,x_5,x_6)=\sum_{N\in\Z^3} v_{N}(x_4)\, e^{2\pi i(n_3x_3+n_5x_5+n_6x_6)}\,,
\end{aligned}
$$
where $N=(n_3,n_5,n_6)$.
Then
\begin{equation}
\label{de5}
\partial_{x_{5}}u=\sum_{N\in\Z^3} 2\pi i n_{5} u_{N}(x_4)\,e^{2\pi i(n_3x_3+n_5x_5+n_6x_6)}\,,
\end{equation}
and the same relation holds for $\partial_{x_{6}}u$,
$\partial_{x_{5}}v$, $\partial_{x_{6}}v$. Hence, by plugging
\eqref{de5} and the other expressions for the derivatives of $u,v$
into equations c., we get
$$
\begin{cases}
&(n_5+x_4n_6)\,u_N(x_4)-n_6\,v_N(x_4)=0\\[3pt]
&n_6u_{N}\,(x_4)+(n_5+x_4n_6)\,v_N(x_4)=0\,,
\end{cases}
$$
for any $N=(n_3,n_5,n_6)\in\Z^3$. If $(n_5+x_4n_6)^2+n_6^2\neq 0$
then $u_N(x_4)=v_N(x_4)=0$. Therefore if $f$ satisfies equations
a. and c. then $f=f(x_3,x_4)$. In particular $f$ must be
$\Z^2$-periodic. By equations b. we immediately get $f\equiv 0$.
Hence $J$ is not admissible. }
\end{ex}
\section{Moduli spaces of admissible almost complex structures}
Let $(M,\omega)$ be a symplectic manifold with vanishing first
Chern class. By using notation of section 2 let
$$
\mathcal{C}_\omega(M)=\{J\in\mbox{End}(TM)\,|\,J^2=-I\,,\,\omega(\cdot,\cdot)=\omega(J\cdot,J\cdot)\,,\,\,\omega(\cdot,J\cdot)>0\}
$$
be the space of $\omega$-calibrated almost complex structures on $M$.\\
Let $J$ be an $\omega$-calibrated almost complex structure; then we
say that $\widetilde{J}\in\mathcal{C}_{\omega}(M)$ is \emph{close}
to $J$ if det$(I-\wJ J)\neq 0$. It is known that the space of
$\omega$-calibrated almost complex structures close to a fixed  $J$
is parametrized by the tangent bundle symmetric endomorphisms
anti-commuting with $J$ and having norm less than 1: namely $\wJ$ is
close to $J$ if and only if there exists a unique $L\in
\mbox{End}(TM)$ such that
$$
\begin{cases}
\begin{aligned}
&\wJ=RJR^{-1}\,,\\
&^tL=L\,,\\
&LJ=-JL\,,\\
&||L||<1\,,
\end{aligned}
\end{cases}
$$
where $R=I+L$ and the norm $||\cdot||$ is taken with respect to
$g_J$ (see \cite{AL}).

Denote by
$$
\mathcal{AC}_{\omega}(M)=\{J\in \mathcal{C}_{\omega}(M)\,|\,\,J
\mbox{ is } \omega-\mbox{admissible}\}\,.
$$
Then the {\em symplectic group}
$$
\hbox{\rm Sp}_{\omega}(M)=\left\{\phi\in\hbox{\rm
Diff}(M)\,\,\vert\,\,\phi^*\omega=\omega \ \right\}
$$
acts on $\mathcal{AC}_{\omega}(M)$ by conjugation
$$
(\phi,J)\longmapsto \phi_*^{-1}\circ J\circ\phi_*\,.
$$
Let
$$
\mathfrak{M}(\mathcal{AC}_\omega(M))=\mathcal{AC}_{\omega}(M)/\mbox{Sp}_\omega(M)
$$
be the relative {\em Moduli Space}.
\subsection{Deformation of $\op_J$}In order to give
a description of $\mathfrak{M}(\mathcal{AC}_\omega(M))$, we have to
describe the behavior of the $\op$ operator for an almost complex
structure $\widetilde{J}$ close to a fixed $J$. We start considering
the following
\begin{prop}
Let $R=I+L$ be an arbitrary isomorphism of $TM$. Then
\begin{equation}
\label{deba}
RdR^{-1}\gamma=d\gamma+[\tau_L,d]\gamma+\sigma_L\gamma
\end{equation}
for any exterior form $\alpha$ in $M$ of positive degree, where:
\begin{itemize}
\item $\tau_L$ is the
zero order derivation defined on the $r$-forms by
$$
\begin{aligned}
\tau_L\gamma(X_1,X_2\dots,X_r)=&\gamma(LX_1,X_2,\dots,X_n)+\gamma(X_1,LX_2,\dots,X_n)+\\[10pt]
&\dots+\gamma(X_1,X_2,\dots,LX_n)\,;
\end{aligned}
$$
\item $[\tau_L,d]=\tau_Ld-d\tau_L$;
\item $\sigma_L$ is the operator defined on the $1$-form as
$$
\sigma_{L}\alpha(X,Y):=\alpha(R^{-1}(N_{L}(X,Y)))\,,
$$
being $N_{L}(X,Y):=[LX,LY]-L[LX,Y]-L[X,LY]+L^2[X,Y]$, and extended on the forms of arbitrary degree by the Leibniz rule.
\end{itemize}
\end{prop}
\begin{proof}
Let $\alpha\in\Lambda^1(M)$ and let $X,Y\in TM$. We have
$$
\begin{aligned}
RdR^{-1}\alpha(X,Y)=&dR^{-1}\alpha(RX,RY)=RX\alpha(Y)-RY\alpha(X)+\alpha(R^{-1}[RX,RY])\\
                   =&X\alpha(Y)-Y\alpha(X)+\alpha([X,Y])+
                   LX\alpha(Y)-LY\alpha(X)\\
                   &+\alpha(R^{-1}[RX,RY]-[X,Y])\\
                   =&d\alpha(X,Y)+ LX\alpha(Y)-LY\alpha(X)
                   +\alpha(R^{-1}[RX,RY]-[X,Y])\,.
\end{aligned}
$$
Moreover
$$
\begin{aligned}
\tau_Ld\alpha(X,Y)=&d\alpha(LX,Y)+d\alpha(X,LY)\\
                  =&LX\alpha(Y)-Y\alpha(LX)+\alpha([LX,Y])+X\alpha(LY)-LY\alpha(X)+\\
                  &\alpha([X,LY])
\end{aligned}
$$
and
$$
d\tau_L\alpha(X,Y)=X\alpha(LY)-Y\alpha(LX)+\alpha(L[X,Y])\,.
$$
Therefore we obtain
$$
\begin{aligned}
(RdR^{-1}-[\tau_L,d])\alpha(X,Y)=&d\alpha(X,Y)+\\
                                 &\alpha(R^{-1}[RX,RY]-[LX,Y]-[X,LY]+L[X,Y]-[X,Y])
\end{aligned}
$$
Now we have
$$
\begin{aligned}
&R(R^{-1}[RX,RY]-[LX,Y]-[X,LY]+L[X,Y]-[X,Y])=\\
&=[RX,RY]-R[LX,Y]-R[X,LY]+RL[X,Y]-R[X,Y]=\\
&=[LX,LY]+[LX,Y]+[X,LY]+[X,Y]-[LX,Y]-L[LX,Y]\\
&\,\,\quad-[X,LY]-L[X,LY]+L[X,Y]+L^2[X,Y]-[X,Y]=\\
&=N_{L}(X,Y)\,,
\end{aligned}
$$
i.e.
$$
R^{-1}[RX,RY]-[LX,Y]-[X,LY]+L[X,Y]-[X,Y])=R^{-1}(N_L(X,Y))
$$
Hence we have
$$
(RdR^{-1}-[\tau_L,d])\alpha(X,Y)=d\alpha(X,Y)+\alpha(R^{-1}(N_L(X,Y)))
$$
which proves the proposition when $\alpha$ is a $1$-form. Since the operators on the two sides of formula \eqref{deba}
satisfy Leibnitz rule, the proof is complete.
\end{proof}
Now we are ready to give the following
\begin{prop}
\label{santiprop} Let $J,\wJ$ be close almost complex structures in
$\mathcal{C}_{\omega}(M)$. Let $\ov{\partial}_J$,
$\ov{\partial}_{\wJ}$ be the $\ov{\partial}$-operators with respect
to $J$, $\wJ$ respectively. Then
\begin{enumerate}
\item[1.] $R\op_{\wJ}f=\op_J f+L\partial_J f\,,$
\vspace{0.2 cm}
\item[2.] $R\op_{\wJ}R^{-1}
  \alpha=\op_{J}\alpha+[\tau_L,d]^{p,q+1}\alpha+\sigma_L^{p,q+1}\alpha\,,$
\end{enumerate}
where $f\in C^{\infty}(M,\C)$, $\alpha\in\Lambda_{\wJ}^{p,q}(M)$,
$\wJ=RJR^{-1}$, $R=I+L$,  $[\tau_L,d]^{p,q+1},\,(\sigma_L)^{p,q+1}$ denote the projection of the bracket
$[\tau_L,d]=\tau_Ld-d\tau_L$ and of the operator $\sigma_L$  on the
space $\Lambda^{p,q+1}_J(M)$, respectively.
\end{prop}
\begin{proof}
1. Let $f\in C^{\infty}(M,\C)$. We have
$$
R\op_{\wJ}f=R(df)^{\widetilde{0,1}}=(Rdf)^{0,1}=(df+Ldf)^{0,1}=\op_Jf+Ldf^{0,1}=\op_Jf+L\partial_Jf^{0,1}\,,
$$
where the subscript $\widetilde{0,1}$ denotes the projection onto $\Lambda^{0,1}_\wJ(M)$.
\vskip0.4cm
2. Let $\alpha\in\Lambda^{p,q}_{J}(M)$. Then we have
$$
\begin{aligned}
R\op_{\wJ}R^{-1}\alpha&=R(dR^{-1}\alpha)^{\widetilde{p,q+1}}=(RdR^{-1}\alpha)^{p,q+1}\\
                      &=(d\alpha)^{p,q+1}+[\tau_L,d]^{p,q+1}\alpha+\sigma_L^{q,p+1}\alpha\\
                      &=\op_J\alpha+[\tau_L,d]^{p,q+1}\alpha+\sigma_L^{q,p+1}\alpha\,,
\end{aligned}
$$
where the subscript $\widetilde{p,q+1}$ denotes the projection onto $\Lambda^{p,q+1}_\wJ(M)$.
\end{proof}
\subsection{The deformation form}\label{santideformationform}
In this subsection we introduce a $(0,1)$-form which will be a useful tool to study infinitesimal
deformations of admissible almost complex structures.

Let $(M,\omega)$ be a symplectic manifold, $J\in\mathcal{C}_\omega(M)$ and $\psi\in\Lambda^{n,0}_J(M)$ nowhere vanishing.
For any endomorphism $L$ of $TM$ anticommuting with $J$ there exist unique forms
$\mu_L(\psi),\,\tau_L(\psi)\in\Lambda^{0,1}(M)$
satisfying the following relations
$$
(\tau_L\ov{A}_J\,\psi)^{n,1}=\mu_L(\psi)\wedge\psi\,,\quad \partial_J\tau_L\,\psi=\gamma_L(\psi)\wedge\psi\,.
$$
\begin{definition}
The $(0,1)$-form
$$
\theta_L(\psi):=\mu_L(\psi)-\gamma_L(\psi)
$$
is called the \emph{deformation form} of $L$.
\end{definition}
\noindent Note that if $J$ is integrable, then $\theta_L(\psi)=-\gamma_L(\psi)$, since $\ov{A}_J=0$.

The following lemma  gives the behavior of  $\theta_L(\psi)$
when the complex volume form $\psi$ changes:
\begin{lemma}
\label{santilemma} Let $\psi,\,\psi'\in\Lambda^{n,0}_J(M)$ be
$\op_J$-closed. Let $\{Z_1,\dots,Z_n\}$ be a local $(1,0)$-frame
and
$\{\zeta_1,\dots,\zeta_n\}$ be the dual frame.\\
Then
\begin{enumerate}
\item[1.]$\mu_L(\psi')=\mu_L(\psi)\,,$ \vspace{0.2 cm}
\item[2.]$\gamma_L(\psi')=\gamma_L(\psi)+\eta(f)\,,$
\end{enumerate}
where $\psi'=f\psi$ and $\eta(f)$ is the $(0,1)$-form defined
locally as
$$
\eta(f)=-\frac{1}{f}\sum_{k,r=1}^n Z_k(f)L_{k\ov{r}}\,\ov{\zeta}_{r}\,,
$$
where
$$
L(Z_{i})=\sum_{k=1}^n L_{\ov{k}i}\,\ov{Z}_{k}\,.
$$
\end{lemma}
\begin{proof}
By definition we have
$$
\begin{aligned}
\mu_L(\psi')\wedge\psi'=&(\tau_L\ov{A}_J\psi')^{n,1}=(\tau_L\ov{A}_Jf\psi)^{n,1}
=f(\tau_L\ov{A}_J\psi)^{n,1}\\
=&f\,\mu_L(\psi)\wedge\psi=\mu_L(\psi)\wedge\psi'\,.
\end{aligned}
$$
Therefore 1. is proved.\\
We have
$$
\begin{aligned}
\gamma_L(\psi')\wedge\psi'&=\partial_J\tau_L\psi'=\partial_J\tau_Lf\psi\\
                                 &=\partial_Jf\wedge\tau_L\psi+f\,\partial_J\tau_L\psi\\
                                 &=\partial_Jf\wedge\tau_L\psi+f\,\gamma_L(\psi)\wedge\psi\\
                                 &=\partial_Jf\wedge\tau_L\psi+\gamma_L(\psi)\wedge\psi'\,.
\end{aligned}
$$
Now we express $\partial_Jf\wedge\tau_L\psi$ in terms of $\psi'$.
With respect to the local (1,0)-frame $\{Z_1,\dots,Z_n\}$ we get
$$
\psi=h\,\zeta_1\wedge\dots\wedge \zeta_n\,,
$$
where $h$ is a local nowhere vanishing smooth function
and
$$
\partial_Jf=\sum_{k=1}^nZ_k(f)\,\zeta_k\,.
$$
Now we have
$$
\begin{aligned}
\tau_L\psi&=h\,L(\zeta_1)\wedge\dots\wedge\zeta_n+\dots+h\,\zeta_1\wedge\dots\wedge
              L(\zeta_n)\\
              &=h\sum_{r=1}^n\{L_{1\ov{r}}\,\ov{\zeta}_{r}\wedge\dots\wedge\zeta_n\}+\dots+
                h\sum_{r=1}^n\{(-1)^{n-1} L_{n\ov{r}}\,\ov{\zeta}_{r}\wedge\dots\wedge\zeta_{n-1}\}\,.
\end{aligned}
$$
Therefore
$$
\begin{aligned}
\partial_Jf\wedge\tau_L\psi=&-\sum_{k,r=1}^nZ_k(f)L_{k\ov{r}}\,\ov{\zeta}_{r}\wedge\psi\\
                               =&-\frac{1}{f}\sum_{k,r=1}^nZ_k(f)L_{k\ov{r}}\,\ov{\zeta}_{r}\wedge\psi'\,.
\end{aligned}
$$
Hence
$$
\gamma_L(\psi')=\gamma_L(\psi)+\eta(f)\,,
$$
i.e. 2. is proved.
\end{proof}
\begin{cor}\label{santicorollario}
Assume that $(M,\omega)$ is compact. Then the deformation form $\theta_L$ does not depend on the choice of the
$\op_J$-closed complex volume form on $M$.
\end{cor}
\begin{proof}
Since $M$ is compact if $\psi,\psi'$ are two complex volume form on $M$ satisfying $\op_J\psi=\op_J\psi'=0$, then
$\psi=c\,\psi'$ for some constant $c$ on $M$. Hence by lemma \ref{santilemma}
$$
\theta_L(\psi)=\theta_L(\psi')
$$
for any $L\in\mbox{End}(TM)$ anticommuting with $J$.
\end{proof}
\subsection{Infinitesimal deformations of admissible complex structures}
In this subsection we compute the infinitesimal deformations of admissible almost complex structures
on a symplectic manifold. We start with the following
\begin{theorem}\label{lemmacentrale}
Let $(M,\omega)$ be a symplectic manifold. Fix $J
\in\mathcal{AC}_\omega(M)$ and consider a smooth curve $J_t$ in
$\mathcal{AC}_\omega(M)$ of almost complex structures close to $J$
satisfying $J_0=0$. Then the derivative $\dot{J}_0$ of $J_t$ at $0$
is given by
$$
\dot{J}_0=2\,JL,
$$
where
$L\in\mbox{\emph{End}}(TM)$ satisfies  the following conditions
$$
L=\,^t\!L\,,\quad LJ=-JL
$$
and for any nowhere vanishing $\op_J$-closed
$\psi\in\Lambda^{n,0}_J(M)$ the $(0,1)$-form $\theta_L(\psi)$ is $\op_J$-exact.

\end{theorem}
\begin{proof}
Let $\psi\in\Lambda^{n,0}_J(M)$ be a nowhere vanishing complex form
satisfying $\op_J\psi=0$.\\ Fix an $\omega$-calibrated almost
complex structure $\wJ$ close to $J$. Then
$$
\wJ=RJR^{-1}\,,
$$
where
$$
R=I+L\,,\quad L=^t\!\!L\,,\quad LJ+JL=0\,,\quad \|L\|<1\,.
$$
The form $R^{-1}\psi\in \Lambda^{n,0}_{\wJ}(M)$ is  nowhere vanishing. Any other section $\psi'$ which
trivializes $\Lambda^{n,0}_{\wJ}(M)$ is a multiple of $\psi$,
namely
$$
\psi'=f\,R^{-1}\psi\,,
$$
with $f\in C^{\infty}(M,\C)$, $f(p)\neq 0$ for
any $p\in M$.\\
Hence the almost complex structure $\wJ$ is $\omega$-admissible if and only if there exists $f\in
C^{\infty}(M,\C)$ such that
$$
\op_{\wJ}(fR^{-1}\psi)=0\,,
$$
where $f\neq 0$.\\
By formulae of proposition \ref{santiprop} we have
$$
\begin{aligned}
R\op_{\wJ}(fR^{-1}\psi)=&R(\op_{\wJ}f\wedge
R^{-1}\psi+f\op_{\wJ}R^{-1}\psi)\\[3 pt]
=&R(\op_{\wJ}f)\wedge\psi+fR\op_{\wJ}R^{-1}\psi\\[3 pt]
=&\op_{J}f\wedge\psi+L\partial_{J}f\wedge\psi+f(\op_{J}\psi+[\tau_L,d]^{n,1}\psi
+(\sigma_L)^{n,1}\psi)\\[3 pt]
=&\op_{J}f\wedge\psi+L\partial_{J}f\wedge\psi+f([\tau_L,d]^{n,1}\psi
+(\sigma_L)^{n,1}\psi)\,,
\end{aligned}
$$
i.e.
\begin{equation}
\label{santiformula}
R\op_{\wJ}(fR^{-1}\psi)=\op_{J}f\wedge\psi+L\partial_{J}f\wedge\psi+f[\tau_L,d]^{n,1}\psi
+f(\sigma_L)^{n,1}\psi\,.
\end{equation}

Let us consider a smooth curve of $\omega$-admissible complex
structures $J_t$ close to $J$, such that $J_0=J$. For any $t$ there
exists $L_t\in \mbox{End}(TM)$ such that if $R_t=I+L_t$, then
$$
J_t=R_tJR_t^{-1}\,,
$$
for $L_tJ+JL_t=0$, $\|L_t\|<1$ . We may assume $L_0=0$. We set $\dot{L}_0=L$. Then
$$
\dot{J}_0=JL-LJ=2\,JL\,.
$$
For any $t$ there exists $f_t\colon
M\to\C$, $f_t\neq 0$, such that
$$
\op_{J_t}(f_tR_t^{-1}\psi)=0\,.
$$
Hence by formula (\ref{santiformula}) $J_{t}$ is
$\omega$-admissible if and only exists $f_{t}\colon M\to \C$, such
that $f_t\neq 0$ and
\begin{equation}
\label{santi_t}
\op_{J}f_t\wedge\psi+L_t\partial_{J}f_t\wedge\psi+f_t
[\tau_{L_t},d]^{n,1}\psi +f_t(\sigma_{L_t})^{n,1}\psi=0\,.
\end{equation}
We may assume without loss of generality  that
$$
\begin{aligned}
&f_0=1\,.
\end{aligned}
$$
Since $\frac{d}{dt}\sigma_{L_t}|_{t=0}=0$, by taking the derivative of  (\ref{santi_t}) at $t=0$ we get
\begin{equation}
\label{santiderivata}
\op_{J}\dot{f}_0\wedge\psi+[\tau_L,d]^{n,1}\psi=0\,.
\end{equation}
Let us compute $[\tau_L,d]^{n,1}\psi$. Since
$d=A_{J}+\partial_J+\op_J+\ov{A}_J$, we have
$$
\begin{aligned}
& (\tau_Ld\psi)^{n,1}=(\tau_L\ov{A}_J\psi)^{n,1}\,,\\
& (d\tau_L\psi)^{n,1}=\partial_J\tau_L\psi\,.
\end{aligned}
$$
By the definition of $\mu_L(\psi)$, $\gamma_L(\psi)$ we have
\begin{equation*}
(\tau_L\ov{A}_J\psi)^{n,1}=\mu_L(\psi)\wedge\psi\,,\quad
\partial_J\tau_L\psi=\gamma_L(\psi)\wedge\psi\,.
\end{equation*}
Hence (\ref{santiderivata}) reduces to
\begin{equation*}
\op_{J}\dot{f}_0\wedge\psi+\mu_L(\psi)\wedge\psi-\gamma_L(\psi)\wedge\psi=0\,,
\end{equation*}
so that \eqref{santiderivata} is equivalent to
\begin{equation*}
\label{santiscalare}
\op_{J}\dot{f}_0+\mu_L(\psi)-\gamma_L(\psi)=0\,,
\end{equation*}
i.e. $\theta_L(\psi)$ is $\op_J$-exact.
\end{proof}
We give the following
\begin{definition}
Let $(M,\omega)$ be a symplectic manifold and $J\in\mathcal{AC}_\omega(M)$. The tangent space to the moduli
space $\mathfrak{M}(\mathcal{AC}_\omega(M))$ at $[J]$ is the vector space
$$
T_{[J]}\mathfrak{M}(\mathcal{AC}_\omega(M)):=\frac{\{\dot{\gamma}(0)\,\,|\,\,\gamma\colon(-\epsilon,\epsilon)
\to\mathcal{AC}_\omega(M)\,\,\mbox{ is a smooth curve s.t.}\,\,\gamma(0)=J\}}{T_J\mathcal{O}_J(M)}\,,
$$
where $\mathcal{O}_J(M)$ denotes the orbit of $J$ under the action of $\mbox{\emph{Sp}}_\omega(M)$.
\end{definition}
As a corollary of theorem \ref{lemmacentrale} we have the following
\begin{prop}\label{centrale}
Let $(M,\omega)$ be a compact symplectic manifold of dimension $2n$
and  let $J\in\mathcal{AC}_\omega(M)$. Then the tangent space to
the Moduli space $\mathfrak{M}(\mathcal{AC}_{\omega}(M))$ at the
point $[J]$  satisfies
$$
T_{[J]}\mathfrak{M}(\mathcal{AC}_\omega(M))\subseteq \frac{\{JL\in
\mbox{\emph{End}}(TM)\,\,|\,\,L=^t\!\!L\,,LJ=-JL\,,
\theta_L\mbox{ is }\op_J-\mbox{exact}
\}}{\{\mathcal{L}_{X}J\,\,|\,\,X\in TM \mbox{ and
}\mathcal{L}_{X}\omega=0 \}}\,.
$$
\end{prop}
\begin{proof}
By corollary \ref{santicorollario}
the $(0,1)$-form $\theta_L$ does not
depend on the choice of the volume form $\psi$. Therefore, by
theorem \ref{lemmacentrale}, if $J_t$ is a smooth curve in $\mathcal{AC}_\omega(M)$ satisfying
$J_0=J$, then $\dot{J}_0=2\,JL$, where $L$ is a symmetric endomorphism of $TM$ anticommuting with $J$ and such that
$\theta_L$ is $\op_J$-exact.
A standard argument shows that
$$
T_J\mathcal{O}_J(M)=\{\mathcal{L}_XJ\,\,|\,\,X\in TM\mbox{ and
}\mathcal{L}_X\omega=0\}
$$
and this completes the proof
\end{proof}
In analogy to the classical case, it is natural to consider the following
\begin{definition}\label{rigido}
The vector space
$$
TV_{[J]}\mathfrak{M}(\mathcal{AC}_\omega(M)):=\frac{\{JL\in
\mbox{\emph{End}}(TM)\,\,|\,\,L=^t\!\!L\,,LJ=-JL\,,
\theta_L\mbox{ is }\op_J-\mbox{exact}
\}}{\{\mathcal{L}_{X}J\,\,|\,\,X\in TM \mbox{ and
}\mathcal{L}_{X}\omega=0 \}}
$$
will be called the \emph{virtual tangent space} to the moduli space $\mathfrak{M}(\mathcal{AC}_\omega(M))$ at $[J]$.\\

An $\omega$-admissible almost complex structure will be called
\begin{itemize}
\item \emph{non-obstructed} if $T_{[J]}\mathfrak{M}(\mathcal{AC}_\omega(M))=TV_{[J]}\mathfrak{M}(\mathcal{AC}_\omega(M))$\,,
\vskip0.2cm
\item \emph{rigid} if $TV_{[J]}\mathfrak{M}(\mathcal{AC}_\omega(M))=\{0\}$.
\end{itemize}
\end{definition}
\vskip0.3cm
We end this section with the following
\begin{rem}
\emph{We observe the following:}
\begin{enumerate}
\item[i)]\emph{In contrast to the case of  deformation of complex structures (see e.g. \cite{K}),  in  this context,
as far as we know,
the fact that $J$ in non-obstructed cannot be interpreted as a
cohomological condition, since the $\op_J$-operator does not give rise to a cohomology on the manifold;}
\vskip0.3cm
\item[ii)]\emph{All the considerations made in this section can be done by considering
 almost symplectic manifolds instead of symplectic manifolds}.
\end{enumerate}
\end{rem}
\section{Admissible complex structures on the Torus}
In this section we apply the results of section 4 to the torus,
computing explicitly the virtual tangent space to
$\mathfrak{M}(\mathcal{AC}_\omega(\T^{2n}))$.

Let $\T^{2n}=\C^n/\Z^{2n}$ be the standard complex torus of real dimension $2n$
and let $\{z_1,\dots,z_n\}$ be coordinates on $\C^n$,
$z_{\alpha}=x_{\alpha}+ix_{\alpha+n}$ for $n=1,\dots, n$.\\
Then
$$
\begin{aligned}
&\omega_n=\frac{i}{2}\sum_{\alpha=1}^n dz_{\alpha}\wedge d\ov{z}_{\alpha}\,,\\
&\psi_n=dz_1\wedge\dots\wedge dz_n
\end{aligned}
$$
define a Calabi-Yau structure on $\T^{2n}$.
Therefore the standard complex structure
$J_n$ is a $\omega_n$-admissible complex structure on $\T^{2n}$.\\
Now we want to deform $J_n$ computing the virtual tangent space
$TV_{J_n}\frak{M}(\mathcal{AC}_{\omega_n}(M))$ to the Moduli space
$\frak{M}(\mathcal{AC}_{\omega_n}(M))$. According to the previous
section, given a symmetric $L\in \mbox{End}(TM)$ that anticommutes
with $J_n$,  we have to write down the (0,1)-form $\gamma_L=-\theta_L$
defined by
$$
\op_J(\tau_L\psi_n)=\gamma_L\wedge\psi_n\,.
$$
Let
$$
L=\sum_{s,r=1}^{n}\{
L_{\ov{r}s}dz_s\otimes\de{}{\ov{z}_{r}}+\ov{L}_{\ov{r}s}\,d\ov{z}_{s}\otimes\de{}{z_r}\}\,,
$$
where $\{L_{s\ov{r}}\}$ are $\Z^{2n}$-periodic
functions. Then we get
$$
\begin{aligned}
\tau_L\psi_n&=L(dz_1)\wedge\dots\wedge
dz_n+\dots+dz_1\wedge\dots\wedge
              L(dz_n)\\
              &=\sum_{r=1}^n\{L_{1\ov{r}}\,d\ov{z}_{r}\wedge\dots\wedge dz_n\}+\dots+
                \sum_{r=1}^n\{
                (-1)^{n-1} L_{n\ov{r}}\, d\ov{z}_{r}\wedge\dots\wedge dz_{n-1}\}\\
              &=\sum_{r,s=1}^n(-1)^{r+1}L_{s\ov{r}}\,d\ov{z}_{r}\wedge
dz_1\wedge\dots\wedge\widehat{dz_s}\wedge\dots\wedge dz_n\,,
\end{aligned}
$$
where $\,\,\widehat{}\,\,$ means that the corresponding term is
omitted. Therefore we obtain
$$
\partial_{J}(\tau_{L}\psi)=-\sum_{r,s=1}^n \de{}{z_s}L_{s\ov{r}}\,d\ov{z}_{r}\wedge\psi\,,
$$
i.e.
$$
\gamma_L=-\sum_{r,s=1}^n\de{}{z_s}L_{s\ov{r}}\,d\ov{z}_{r}\,.
$$
Then if $J_t\in\mathcal{AC}_\omega(M)$ is a smooth curve satisfying $J_0=J_n$, then $\dot{J}_0=2\,J_nL$, where
$$
L=^{t}\!\!L\,,\quad J_{n}L+LJ_{n}=0
$$
and the $(0,1)$-form
$$
\gamma_L=-\sum_{r,s=1}^n\de{}{z_s}L_{s\ov{r}}\,d\ov{z}_{r}
$$
is $\op_{J_n}$-exact.
In order to compute $TV_{J_n}\frak{M}(\mathcal{AC}_\omega(\T^{2n}))$,
we have to write down the Lie derivative $\mathcal{L}_{X}J_{n}$,
for $X\in\mbox{End}(TM)$, such that
$\mathcal{L}_X\omega_n=0$. Let
$X=\sum_{r=1}^{2n} a_r \de{}{x_{r}}$ be a real vector field on
$\T^{2n}$, then a direct computation gives
$$
\mathcal{L}_{X}\omega_n=0\iff
\begin{cases}
\begin{array}{lll}
\displaystyle\de{}{x_r}a_{n+s}-\de{}{x_s}a_{n+r}&=&0\\[9pt]
\displaystyle\de{}{x_r}a_{s}+\de{}{x_{n+s}}a_{n+r}&=&0\\[9pt]
\displaystyle\de{}{x_r}a_{n+s}-\de{}{x_s}a_{n+r}&=&0
\end{array}
\end{cases}
$$
for $r,s=1,\dots,n$. Furthermore
$$
\begin{aligned}
\mathcal{L}_{X}(J_{n})\Big(\de{}{z_{s}}\Big)=
&\sum_{r=1}^{2n}-i\de{a_{r}}{z_{s}}\de{}{x_{r}}+\de{a_{r}}{z_{s}}J_n\Big(\de{}{x_{r}}\Big)\\
                            =&\sum_{r=1}^{2n}-i\de{a_r}{z_s}\Big(\de{}{x_r}+iJ_n\de{}{x_r}
                            \Big)\\
                            =&\sum_{r=1}^{n}-i\de{a_r}{z_s}\Big(\de{}{x_r}-i\de{}{x_{r+n}}\Big)-
                             i\de{a_{r+n}}{z_{s}}\Big(\de{}{x_{r+n}}-i\de{}{x_{r}}\Big)\\
                            =&-\sum_{r=1}^{n}\de{}{z_s}(ia_{r}+a_{r+n})\Big(\de{}{x_r}+i\de{}{x_{r+n}}\Big)\\
                            =&-2\sum_{r=1}^{n}\de{}{z_s}(ia_{r}+a_{r+n})\de{}{\ov{z}_{r}}\,,
\end{aligned}
$$
i.e.
$$
\mathcal{L}_{X}(J_{n})=-2\sum_{r,s=1}^{n}\Big\{\de{}{z_s}(ia_{r}+a_{r+n})\,dz_{s}\otimes\de{}{
\ov{z}_{r}}
+\de{}{\ov{z}_s}(-ia_{r}+a_{r+n})\,d\ov{z}_{s}\otimes\de{}{z_{r}}\Big\}\,.
$$
Therefore $L=\mathcal{L}_{X}J_{n}$, with $\mathcal{L}_X\omega=0$ if and only if
\begin{equation}
\label{imdelta}
L_{\ov{r}s}=2\de{}{z_s}(a_{r}-ia_{r+n})\,,
\end{equation}
where $a_r\,,a_{r+n}$ are periodic functions on  $\R^{2n}$ satisfying
$$
\begin{cases}
\begin{array}{lll}
\displaystyle\de{}{x_r}a_{n+s}-\de{}{x_s}a_{n+r}&=&0\\[9pt]
\displaystyle\de{}{x_r}a_{s}+\de{}{x_{n+s}}a_{n+r}&=&0\\[9pt]
\displaystyle\de{}{x_r}a_{n+s}-\de{}{x_s}a_{n+r}&=&0\,.
\end{array}
\end{cases}
$$

Now we show that the standard complex structure $J_n$ on $T^{2n}$ is not rigid.\\
In order to see this, we show that if $L$ belongs to tangent space to the orbit of $J_n$ at $J_n$ and it has
constant coefficients, then $L$ is zero.\\
Let $J_nL\in \mbox{End}(T\T^{2n})$, where $L$ is symmetric, anticommutes with $J_n$
and it is  such that
$L_{\ov{r}s}$ are constant functions.\\
By equation \eqref{imdelta} there exists $X\in TM$ such that
$J_{n}L=\mathcal{L}_X(J_n)$
if and only if
$$
\begin{aligned}
(J_{n}L)_{\ov{r}s}&=2\de{}{z_s}(a_{r}-ia_{r+n})=\de{}{x_s}(a_{r}-ia_{r+n})-i\de{}{x_{s+n}}(a_{r}-ia_{r+n})\\
           &=\Big(\de{a_r}{x_s}-\de{a_{r+n}}{x_{s+n}}\Big)-i\Big(\de{a_{r+n}}{x_s}+\de{a_r}{x_{s+n}}\Big)\,.
\end{aligned}
$$
Therefore
$$
\begin{cases}
\displaystyle\de{a_r}{x_s}-\de{a_{r+n}}{x_{s+n}}&=\mbox{constant}\\[7pt]
\displaystyle\de{a_{r+n}}{x_s}+\de{a_r}{x_{s+n}}&=\mbox{constant}\,,
\end{cases}
$$
that imply
$$
\de{^2a_r}{x_s^2}+\de{^2a_r}{x_{s+n}^2}=0\,
$$
for any $r,s=1,\dots,n$.\\
It follows that the $\{a_r\}$ are harmonic functions on the
standard torus $\T^{2n}$ and then they are constant. Therefore any
constant $0\neq L\in \mbox{End}(T\T^{2n})$ anticommuting with $J_n$
defines a non-trivial element of
$T_{[J_n]}\frak{M}(\mathcal{AC}_{\omega_n}(M))$. Moreover any
constant endomorphisms $L_1$, $L_2$ of such type give rise to
different elements of $T_{[J_n]}\frak{M}(\mathcal{A}\Tk(M))$.
Hence $J_n$ is not rigid.


\begin{thebibliography}{12}
\bibitem{AD}
Apostolov V., Draghici T.: The curvature and the integrability
of almost-K\"ahler manifolds: a survey,
\emph{Symplectic and contact topology: interactions and perspectives} (Toronto, ON/Montreal, QC, 2001),  25-53,
Fields Inst. Commun., 35,
Amer. Math. Soc., Providence, RI, 2003.
\bibitem{AL}
\emph{Holomorphic curves in symplectic geometry.} Edited by Mich\'ele Audin
and
Jacques Lafontaine. Progress in Mathematics, 117. Birkhäuser Verlag, Basel, 1994. xii+328 pp.
\bibitem{CS}
Chiossi S., Salamon S.: The intrinsic torsion of $\rm SU(3)$
 and $G\sb 2$ structures.
Differential geometry, Valencia, 2001, 115--133, World Sci.
Publishing, River Edge, NJ, 2002.
\bibitem{dB2} de Bartolomeis P.: Geometric Structures on Moduli Spaces of Special Lagrangian Submanifolds,
\emph{Ann. di Mat. Pura ed Applicata}, IV, Vol. CLXXIX, (2001),
pp. 361--382.
\bibitem{dBT} de Bartolomeis P., Tomassini A.:
On the Maslov Index of Lagrangian Submanifolds of Generalized
Calabi-Yau Manifolds, to appear in {\em Int. J. of Math.}.
\bibitem{CFG}
Cordero L. A., Fern\'andez M., Gray A.
Symplectic manifolds with no K\"ahler structure,
\emph{Topology} 25 (1986), no. 3, pp. 375--380.
\bibitem{FG}
Fino A., Grantcharov G.: Properties of manifolds with skew-symmetric torsion
and special holonomy,
\emph{Adv. Math.} 189 (2004), no. 2, pp. 439--450.
\bibitem{Gau}
Gauduchon P.: Hermitian connections and Dirac operators.
\emph{Boll. Un. Mat. Ital. B} (7) 11 (1997), no. 2, suppl., pp. 257--288.
\bibitem{G}
Gromov M.: Pseudoholomorphic curves in symplectic manifolds.
\emph{Invent. Math.} 82 (1985), no. 2, pp. 307--347.
\bibitem{GHJ} Gross M., Huybrechts D., Joyce D.:
{\em Calabi-Yau manifolds and related geometries}, Lectures from
the Summer School held in Nordfjordeid, June 2001. Universitext.
Springer-Verlag, Berlin, 2003. viii+239 pp.
\bibitem{H1} Hitchin N.J.: {\em Stable forms and special metrics},
Global differential geometry: the mathematical legacy of Alfred
Gray (Bilbao, 2000), 70--89, Contemp. Math., 288, Amer. Math.
Soc., Providence, RI, 2001.
\bibitem{H}
Hitchin N.: Generalized Calabi-Yau manifolds,
\emph{Q. J. Math.} 54 (2003), no. 3, pp. 281--308.
\bibitem{HE}
Hopf E.: Elementare Bemerkungen ueber die Loesung
parzieller Differentialgleichungen zweiter Ordnung von elliptischen
Typus, \emph{Sitzungber. Preuss. Akad. Wiss. phys. math. Kl}. 19
(1927), pp. 147--152.
\bibitem{J1} Joyce, Dominic D.: {\em Compact manifolds with special
holonomy}, Oxford Mathematical Monographs. Oxford University
Press, Oxford, 2000. xii+436 pp.
\bibitem{K}
Kodaira K.: \emph{Complex manifolds and deformation of complex
structures},
translated from the 1981 Japanese original by Kazuo Akao. Reprint of
the 1986 English edition.
Classics in Mathematics. Springer-Verlag, Berlin, 2005. x+465 pp.
\end{thebibliography}
\end{document}